%% file: ChenTongJTSP12subFinal.tex
\documentclass[conference,10pt]{IEEEtran}
\IEEEoverridecommandlockouts
\makeatletter
\def\ps@headings{%
\def\@oddhead{\mbox{}\scriptsize\rightmark \hfil \thepage}%
\def\@evenhead{\scriptsize\thepage \hfil \leftmark\mbox{}}%
\def\@oddfoot{}%
\def\@evenfoot{}}
\makeatother 
\usepackage{url,amsmath,graphicx,amssymb,stfloats,subfigure}
\usepackage{epsfig,epsf,psfrag,amssymb,amsfonts,latexsym,slashbox,graphicx,bm,cite,url}
\usepackage[dvips]{color}
\usepackage{mathrsfs}

\newcommand{\mbbE}{\mathbb{E}}

\def\ie{{\it i.e.,\ \/}}

\newcommand{\Cmsc}{\mathscr{C}}

\newtheorem{theorem}{Theorem}
\newtheorem{lemma}{Lemma}

\begin{document}
\title{Distributed Learning and Multiaccess\\ of On-Off Channels
\thanks{Shiyao Chen and Lang Tong are with the School of Electrical and Computer Engineering,
Cornell University, Ithaca, NY 14853. Email: {\tt
\{sc933,lt35\}@cornell.edu}.\hfill \break The work is supported in
part by the National Science Foundation under award CCF 1018115 and
the Army Research Office  under award W911NF1010419.
}}
\author{Shiyao Chen and Lang Tong}



\maketitle

%
\input Intro2

\section{System Model and Assumptions}\label{sec:model}
\input ChannelModel

\section{A Distributed Learning and Access Policy}\label{sec:scheme}
\input Policy

\section{Main Results}\label{sec:schemeanalysis}
\input Analysis

\section{Numerical Results}

We conduct numerical simulations for various channel and user scenarios.

\subsection{Simulation setup}

We adopt geometric distribution for the channel on and off period lengths. For homogeneous channels situation, the average on and off period lengths
are $\mu^{\mbox{\tiny on}}=3.23$ and $\mu^{\mbox{\tiny off}}=1.43$ (the long term channel available fraction is $\eta=0.693$).
For heterogeneous channels situation, there are half of the channels with channel parameters
$\mu^{\mbox{\tiny on}}=3.23$ and $\mu^{\mbox{\tiny off}}=1.43$ ($\eta=0.693$), and the other half with channel parameters
$\mu^{\mbox{\tiny on}}=3.23$ and $\mu^{\mbox{\tiny off}}=4.3$ ($\eta=0.429$).

The initial detection period length is $L_0=24$ slots, each time incremented by $C=12$ slots. The entire horizon is taken to be 5000 slots, and the number of Monte Carlo runs is 20.

\subsection{Increasing the Number of Users $K$}

\begin{figure}
\centering
\includegraphics[width=2.2in]{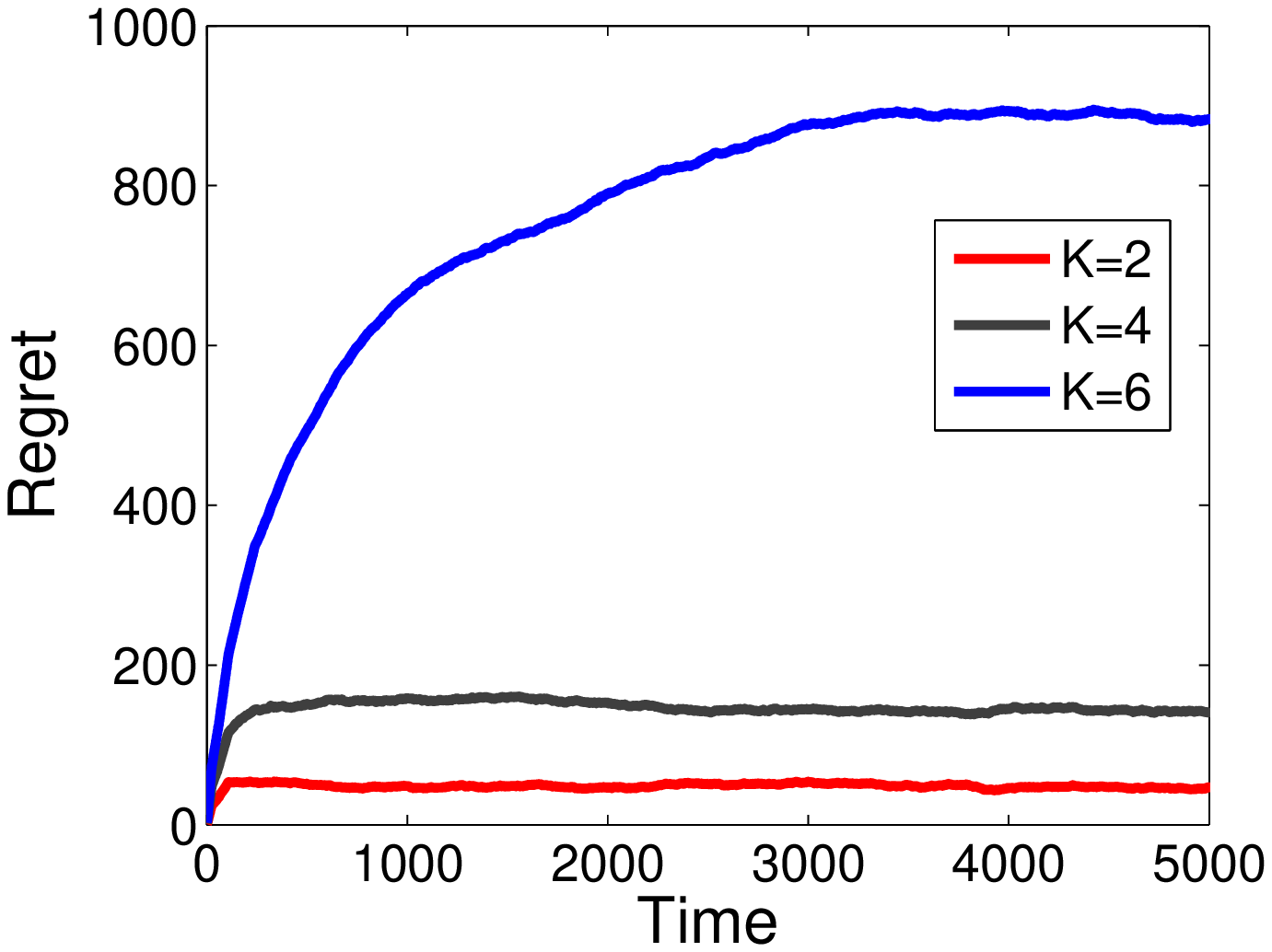}
 \caption{Homogeneous channels. $K=2,4,6$, $N=6$}\label{fig:increasingK}
\end{figure}

We simulate the effect of increasing number of users with $N=6$ and $K=2,4,6$, and the simulated regret is depicted in Fig. \ref{fig:increasingK}.
All the users have targeted individual throughput $r=0.5$.

As predicted in Theorem \ref{thm:regret}, the simulated expected regret indeed levels off eventually, verifying the result of finite expected regret.
The three curves for $K=2,4,6$ in Fig. \ref{fig:increasingK} clearly show an increasing trend of the expected regret and the time it takes for the expected regret to converge when $K$ increases. This trend is quite intuitive; when there are more users, the entire process takes much longer.

\subsection{Fixed vs Increasing Detection Period Length}

\begin{figure}
\centering
\includegraphics[width=2.2in]{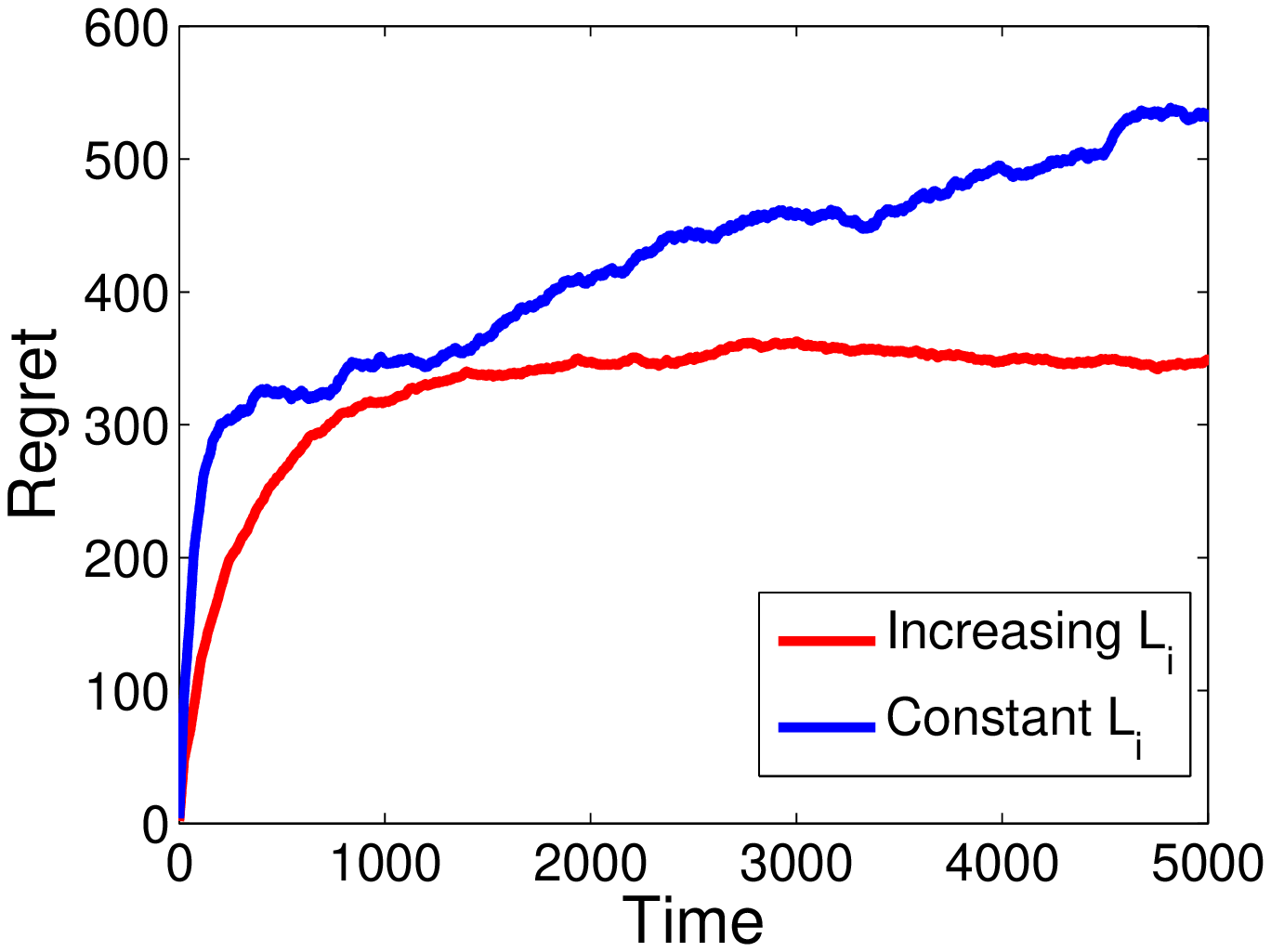}
 \caption{Homogeneous channels. $N=K=4$.}\label{fig:Ho}
\end{figure}

\begin{figure}
\centering
\includegraphics[width=2.2in]{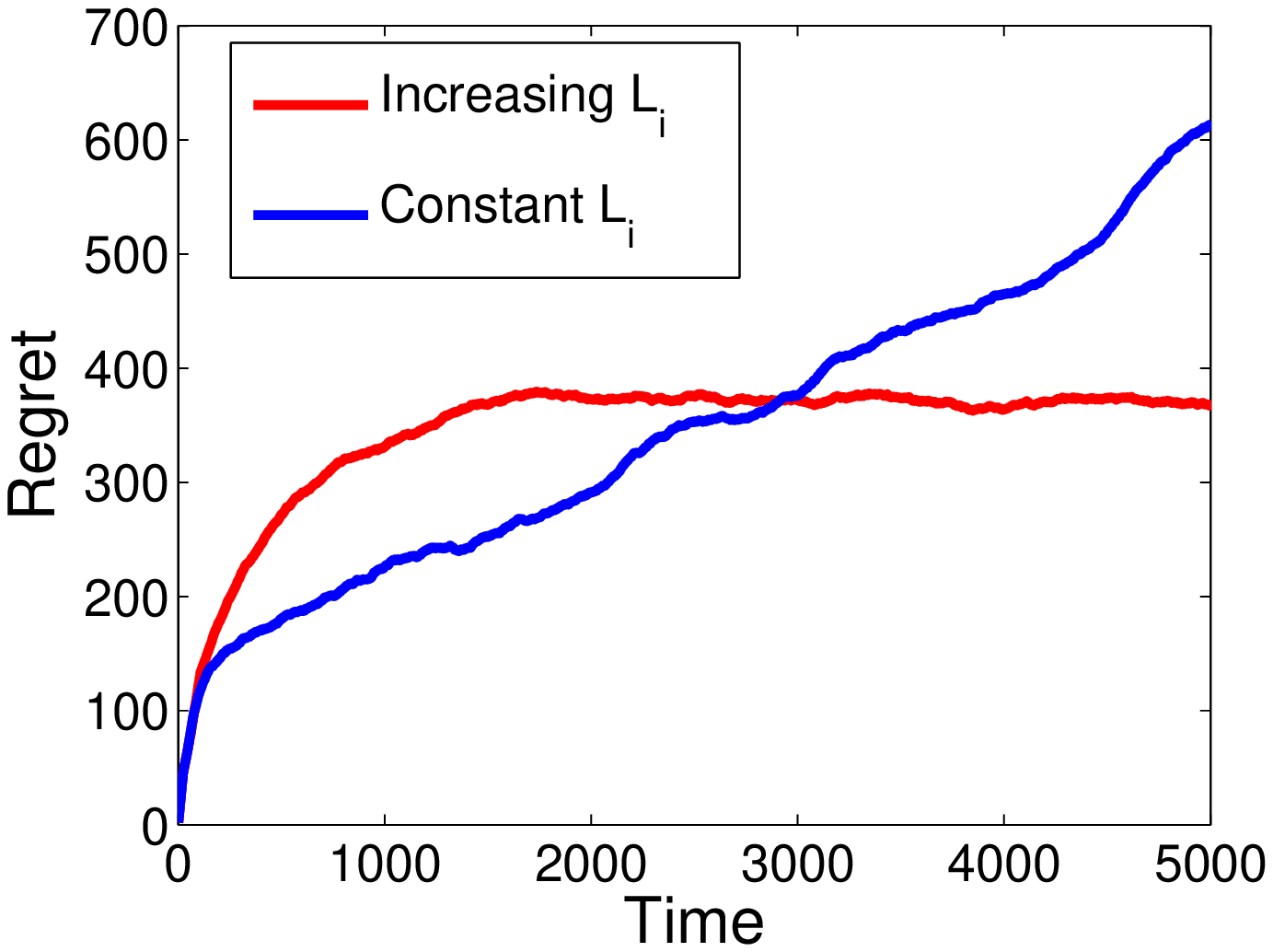}
 \caption{Heterogeneous channels. $N=K=4$.}\label{fig:He}
\end{figure}

To compare the impact of fixed and increasing detection period length, we simulate the situation with initial detection period length 24, and incremental of 12 slots and 0 slot each time.
The simulated expected regret curves are shown in Fig. \ref{fig:Ho} and \ref{fig:He} for $N=K=4$ with fixed and increasing detection period length.

The expected regret associated with increasing detection period length outperforms the counterpart with fixed detection period length. This comparison demonstrates the necessity of the increasing detection period length structure for the desired finite expected regret, and justifies the rationale of establishing the exponential decay of the error probability in detection.

\section{Conclusion}

We consider in this paper the problem of distributed learning and multiaccess of orthogonal channels. We have shown that perfect orthogonalization can be achieved by a distributed and asynchronous learning and access policy in the sense that the throughput region of a centralized scheme with fixed channel allocation can be achieved. In fact, we have established that the expected regret of the proposed distributed scheme with respect to a centralized scheme is finite.
Interesting future directions include analyzing upper bound of the expected regret with system parameters as well as bounds of the expected regret within a finite time horizon.
\section*{Appendix}
\input Appendix

\bibliographystyle{IEEEbib}
{\small\bibliography{refs}}

\end{document}

%% file: Intro2.tex
\begin{abstract}
The problem of distributed  access of a set of $N$ on-off channels
by $K \le N$ users is considered. The channels are slotted and
modeled as independent but not necessarily identical alternating
renewal processes.  Each user decides to either observe  or transmit
at the beginning of every slot. A transmission is successful only if
the channel is at the on state and there is only one user
transmitting. When a user observes, it identifies whether a
transmission would have been successful had it decided to transmit.
A distributed learning and access policy referred to as alternating
sensing and access (ASA) is proposed. It is shown that ASA has
finite expected regret when compared with the optimal centralized
scheme with fixed channel allocation.
 \\[0.2em]
{\em Index terms}---Multiaccess, Distributed learning, Cognitive radio networks.
\end{abstract}

\section{Introduction}\label{sec:intro}
\IEEEPARstart{T}{he} problem considered in this paper, in its more
general form,  is related to distributed allocation of $N$
independent and randomly available resources among $K$ agents. By
distributed allocation we mean that there is no central controller
assigning resources to agents, and each agent acts on its own
without communicating with others.  We are interested in whether
there is a distributed access policy that, through taking actions
and learning from outcomes of the actions, achieves the utilization
of resources comparable to that of the optimal centralized
allocation.

We study the above in the context of multiaccess of $N$ random
on-off channels by $K$ independent users.  We are interested in
whether any distributed learning and access policy is necessarily
penalized by the lack of coordination. The performance measure of
interest is throughput---the fraction of time that transmissions are
successful. For a $K$ user multiaccess system, the throughput is
defined by a vector $\mathbf{r}=(r_1,\cdots, r_K)$ where $r_i$ is
the throughput of user $i$.  If $\mathbf{r}$ can be achieved by a
central controller who assigns a channel to each user, we would like
to achieve the same by letting users act independently on their own.

We should point out at the outset that when there are more users
than the number of channels, \ie $K>N$, the throughput region
achievable by the optimal distributed access is, in general,
strictly smaller than that by the optimal centralized scheme.  This
can be seen from the case when $N=1$, which corresponds to the
classical random access problem. The centralized scheme achieves the
maximum sum-rate of $1$ packet per slot. The  celebrated slotted
ALOHA protocol can be used to achieve, asymptotically as
$K\rightarrow \infty$, the sum rate of $e^{-1}$. Although the
optimal policy of distributed random access for this case is
unknown, it is well known \cite{Gallager:book} that distributed
random access cannot achieve the throughput of the optimal
centralized channel allocation scheme.

Thus we restrict ourselves to the case when $N \ge K$. Here it is no
longer obvious that distributed multiaccess performs strictly worse
than the centralized counterpart. Because a user is not restricted
to transmitting on one particular channel, it can search for
opportunities elsewhere to avoid colliding with others. Intuitively,
as $N$ increases, conflicts among users diminish, and users may be
able to orthogonalize themselves to avoid collision. Even when
$K=N$,  a user can learn where other users are transmitting and act
accordingly to avoid collision.

Beside collisions among competing users, we also consider a specific
nontrivial channel imperfection.  In particular, we assume that the
channels are independent on/off random processes.  Such a random
channel model arises naturally from channel fading in wireless
systems. In the context of multiaccess of cognitive radios
\cite{Zhao&Sadler:SPM}, this model includes the situation when a
channel is unavailable when it is occupied by another user of higher
priority. In both cases, it is difficult for a user to identify whether the
transmission failure is caused by
collision with another user or by channel fading.  Learning in such
an uncertain environment cannot be perfect. It is therefore not
obvious that mistakes in learning only cause negligible performance
loss.

Like many online learning problems in uncertain environments, to achieve the best performance requires careful tradeoffs between exploration and exploitation.  The results presented in this paper is an instance of such tradeoff that balances sensing and transmission.

\subsection{Summary of Results}\label{sec:result}
The detailed system model and assumptions are given in
Section~\ref{sec:model}. Here we outline the context of the problem
and summarize our main results. We consider $N$ independent but not
necessarily identical on-off slotted channels.  Our results apply to
more general settings, but at the moment, it is sufficient to think
these channels as independent Bernoulli channels with probability
$\eta_i$ that the $i$th channel is at the on state.  Let
$\boldsymbol{\eta}=(\eta_1,\ldots,\eta_N)$.

For a $K$ user multiaccess system,
it is obvious that any point in
\begin{equation}
\overline{\mathscr{R}}=\{(r_1,\ldots,r_K)\mid
r_{(i)}\leq\eta_{(i)},1\leq i\leq K\},
\end{equation}
can be achieved by a central controller with fixed channel allocation, where $r_{(i)}$ and $\eta_{(i)}$ are the
ordered list of $r_i$ and $\eta_i$, respectively. Indeed $\overline{\mathscr{R}}$ is the largest
achievable region by a central controller under {\em fixed channel allocation} without time sharing arrangement
 and without using channel state realizations.

The main result of this paper is to show that, under the model specified in Section~\ref{sec:model},
every point in $\overline{\mathscr{R}}$ is achievable by a {\em decentralized} access policy.
This result is established by constructing a distributed learning
and access policy executed independently by all users.  The policy
alternates between sensing and access periods, hence referred to as
the alternating sensing and access (ASA) policy.

The throughput result above is a direct consequence of a more
refined analysis based on the notion of regret between the total
number of successful transmissions up to slot $n$ of the optimal
centralized scheme $\overline{S}_i(n)$ and that of the distributed
scheme $S_i(n)$ proposed here.  We show in Theorem
\ref{thm:regret} that, if $\mathbf{r} \in \overline{\mathscr{R}}$, the
expected regret of ASA approaches to a constant, \ie
\begin{equation} \label{eq:O}
\mbbE(\overline{S}_i(n)-S_i(n))\sim O(1).
\end{equation}

We should point out that when channels are homogeneous, \textit{i.e.}, $\eta_i = \eta$ for all $i$, $\overline{\mathscr{R}}$ is the largest throughput region achievable by a centralized controller, and ASA matches with the optimal centralized allocation.  For heterogeneous channels,  while $\overline{\mathscr{R}}$ is the largest throughput region by a centralized controller under  fixed channel allocation, it is not necessarily convex.  Therefore, $\overline{\mathscr{R}}$ may be enlarged by a central controller through time sharing.

In comparing with centralized access policies, we exclude the possibilities of time sharing arrangement, which is a loss of generality.  Such a loss, however, is inevitable since the optimal distributed access policy cannot in general achieve the convexified throughput region\footnote{This can be seen by considering a two channel case when $\eta_1=1$ and $\eta_2=\epsilon$. As $\epsilon \rightarrow 0$, a centralized access can accommodate a rate vector of $(0.5,0.5)$ by time sharing on a single channel, which cannot be achieved by any decentralized access policy.}.

\subsection{Related Work}

The problem of orthogonalizing multiple coexisting users in a distributed manner through learning and individual actions has been studied as a decentralized learning of multi-armed bandit (MAB) processes involving multiple players in \cite{Liu&Zhao:SP}.  Essentially the same problem has also be studied for the multiaccess problem in multiuser cognitive radio systems \cite{Lai&Jiang&Poor:08Asilomar,Anandkumar&Michael&Tang:10INFOCOM}. Further development of the prioritized and fair access cases is provided in \cite{Gai&Krishnamachari:11GLOBCOM,Dai&Gai&Krishnamachari:12INFOCOM}. 

There are similarities and significant differences between these ``MAB approaches'' and that considered in this paper.
The MAB formulation involves independent random processes, often
assumed independent and identically distributed (iid) in time but
may also be Markovian. Each process is associated with an {\em
unknown} deterministic parameter. Lai and Robbins considered the
single user (non-Bayesian) MAB problem, aiming to maximize the
accumulated reward using knowledge learned from the outcome of past
plays \cite{Lai&Robbins:85APP}.  The problem falls in the category
of ``learning through doing''.

The centralized  multiuser version of the MAB problem was considered
in \cite{Anantharam&Varaiya&Walrand:87TAC} as a single user MAB
problem but with the possibility of simultaneously playing multiple
arms.  The decentralized MAB problem was addressed explicitly in
\cite{Liu&Zhao:SP} and in the context of cognitive radio systems in
\cite{Lai&Jiang&Poor:08Asilomar,Anandkumar&Michael&Tang:10INFOCOM}.
Typically, learning in the MAB problem refers to learning which arms
are more favorable to play. The regret of the order optimal
distributed learning with respect to the oracle player often
increases with the number of plays as $O(\log n)$, unlike that in
Eq. (\ref{eq:O}).

The problem considered in this paper does not belong to the category
of MAB problem although it shares some common characteristics with
the MAB formulation. We highlight here three main differences.
First,  unlike the MAB problems, the parameters of the underlying
random processes are known.   Thus there is little ambiguity on
which channels are favorable for transmissions.  Learning in this
context deals  with searching for appropriate channels to transmit,
not knowing (for certain) the presence of other users.

Second, while the objective of MAB involves maximizing reward, we
are interested in whether a set of rewards can be achieved through
learning and transmission; each user does not try to maximize its
throughput.

Third, the uncertainty associated with the presence of other users
is a  key distinction between the problem treated here and the MAB
formulation.  For the multiuser MAB problem, the presence of other
players are certain whenever two players play the same arm. In our
case, a failed transmission may be caused either by collision or by
the fact that the channel is off.

A related problem is learning parameters of multiple independent processes when a user can choose where and when to observe a particular process \cite{Long&etal,Liang&Liu&Yuan:10ITA,Tehrani&Tong&Zhao:12SP}. Without actively engaging with other users, such formulations are more akin to the classical parameter estimation problems, not one of  ``learning through action'' studied in our and the MAB formulations.
%
%

%% file: ChannelModel.tex
The multiaccess system considered includes $N$ channels, $K$
distributed users, and a basestation. We specify their roles in
their interactions and assumptions made in this paper.

\subsection{Channel Model}
The $N$ channels are slotted and statistically independent.  We consider a slot atomic, which means that it cannot be divided further so that multiple
actions can be taken within one slot.  The channel state of each channel in a slot is either ``on'' or ``off'' with ``on'' indicating that the channel can be used for transmission and ``off'' otherwise.  The state  of each channel  is therefore a discrete-time binary process for which we model it as a {\em  renewal sequence} alternating between consecutive ``on'' and ``off'' periods.

The distribution of the on(off) period of channel $i$ is
$F_i^{\mbox{\tiny on}}$($F_i^{\mbox{\tiny off}}$). We assume that
$\mathbb{E}e^{\theta U}$ is well defined for some $t>0$, where $U$ is
distributed as $F_i^{\mbox{\tiny on}}$($F_i^{\mbox{\tiny off}}$).
Denote the mean of $F_i^{\mbox{\tiny
on}}$($F_i^{\mbox{\tiny off}}$) by $\mu_i^{\mbox{\tiny
on}}$($\mu_i^{\mbox{\tiny off}}$), and the long term
fraction of on periods of channel $i$ by $\eta_i =
\mu_i^{\mbox{\tiny on}}/(\mu_i^{\mbox{\tiny on}} +
\mu_i^{\mbox{\tiny off}})$.

\subsection{User Action and Feedback}
Users act independently and persistently, each aimed at achieving
some fixed throughput target. They do not have a synchronized
starting slot; they may enter the system at different times.

A user makes the decision either to access a channel or to sense a
particular channel at the beginning of slots based on the outcomes
of its own past actions. If the user decides to take the action of
accessing channel $i$ in a slot, it transmits a packet to the
basestation over channel $i$. If the action of observing channel $i$
is taken, it monitors channel $i$ in the slot.

When a user transmits over a particular channel to the basestation,
it receives a {\em binary} feedback at the end of the slot from the
basestation indicating whether the transmission is successful.  When
the transmission is successful, we call the channel over which the
transmission occurred {\em available,} which means that the channel
is at the on state and no other user transmits. The user receives a
feedback symbol ``a''. The transmission fails if the channel is at
the off state or when multiple transmissions occur at the same time.
In this case, we call the channel {\em unavailable,} and the user
receives a feedback symbol ``$\bar{a}$''.  We note that the binary
feedback does not specify which type of failure occurred to the
transmission.

If a user decides to sense a particular channel in a slot, it
observes the channel and decides whether the channel is available,
\ie  whether a transmission would have been successful had the user
decided to transmit. The outcome of the sensing action is again
binary with ``$a$'' indicating that the channel is available  and
``$\bar{a}$'' the opposite.

Note that the information obtained by a user in a slot through
observation is identical to that through the feedback of a
transmission. The difference is that there is a potential reward or
damage caused by transmission.

%% file: Policy.tex
In this section, we present a learning and access policy referred to
as alternating sensing and access (ASA) policy.  We show later in
Section~\ref{sec:schemeanalysis} that ASA achieves the same
throughput region as the optimal centralized scheme with fixed channel allocation.

The process of distributed orthogonalization is dynamic. Two users
collide, which may cause one or both switch to a separate channel.
The switch  may cause further collisions with others.  Because a
user cannot be certain that a failed transmission is caused by
collision, it may decide to switch to a different channel when in
fact the failed transmission is caused by channel.  The key of the
learning and access policy presented here is to mix the actions of
transmission and observation to reduce collisions and
 recover when collisions occur. Here we have a case of dynamic learning where a balance between exploitation and exploration must be made.

\begin{figure}[tb]
\centering\begin{psfrags}
\psfrag{a1}[c]{\footnotesize Channel $i$}
\psfrag{a2}[l]{\footnotesize occupied and tail}
\psfrag{a3}[l]{\footnotesize other user(s)}
\psfrag{a5}[c]{\footnotesize no other users}
\psfrag{a4}[l]{\footnotesize vacant or head}
\psfrag{s}[c]{\footnotesize channel selection}
\psfrag{t}[c]{\footnotesize sensing}
\psfrag{p}[c]{\footnotesize access}
\includegraphics[width=2.3in]{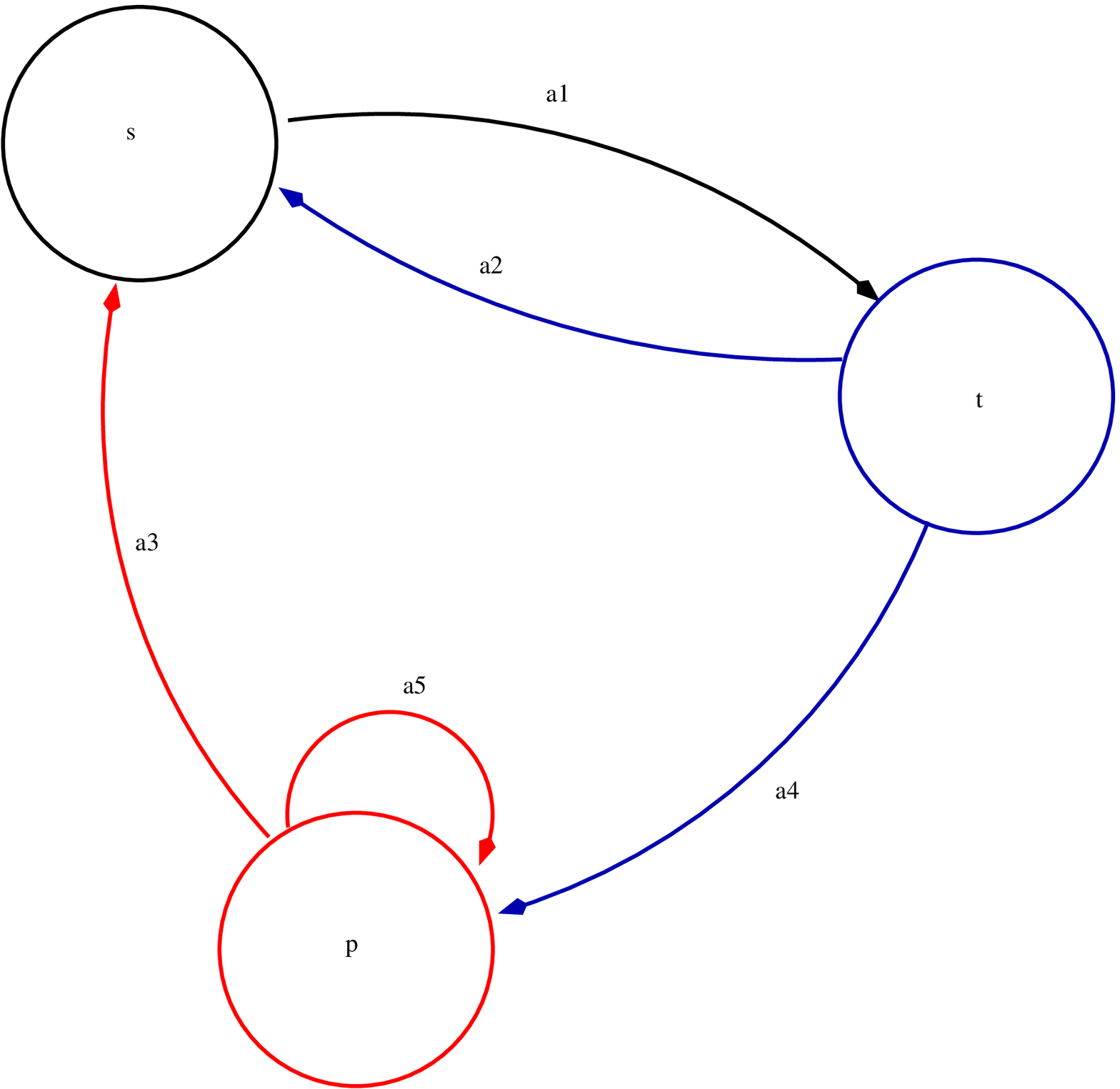}
\end{psfrags}
\caption{State diagram of alternating sensing and access (ASA) policy}\label{fig:diagram}
\end{figure}

\subsection{ASA Policy State and State Transition}
Every user executes the same ASA policy independently. The structure of ASA is illustrated in Fig.~\ref{fig:diagram} where the policy traverses among three policy states: channel selection, sensing, and access. We describe the function of ASA at each state as we follow one user traversing through various states.

We focus on user $i$ who just enters the system, wishing to
communicate at the rate of $r_i$.  User $i$ starts at the channel
selection state knowing that there is a set of channels $\Cmsc_i =
\{k: \eta_k \ge r_i\}$ that can accommodate her rate of
communications. She selects randomly with equal probability one of
the channels as her initial  candidate for access.  With that
choice, she proceeds to  the {\em sensing state}.

At the sensing state, the user senses the channel for a period of
consecutive $L_s$ slots. At the end of the {\em sensing period,} a
detection is performed to test the hypothesis whether the channel is
unoccupied. If the user believes that the channel is unoccupied (she
may be wrong of course), she enters the {\em  access state}. If, on
the other hand, the test result is that the channel is occupied by
another user, she flips a fair coin to further decide whether she
should search for opportunities in other channels, or still enter
the {\em access state} to show her presence to other users. If a
tail shows up, then the user returns to the channel selection state
(as described by ``occupied and tail'' in Fig. \ref{fig:diagram}).
There, again, she chooses randomly another channel from $\Cmsc_i$
(with replacement). Otherwise, if a head shows up, the user still
enters the {\em access state} to transmit and let other users be
aware of her presence (as described by ``vacant or head'' in Fig.
\ref{fig:diagram}).

At the access state, user $i$ transmits in channel $k$ with probability $q_i=r_i/\eta_k$ for a
period of consecutive $L_t$ slots where $q_i$ is chosen to achieve
the desired throughput target $r_i$.  At the end of each slot, user $i$ receives
a feedback. At the end of the {\em transmission period}, a
hypothesis test is made to check if she has been colliding with
another user.  If the test result is that another user is accessing
the channel at the same time (again she may be wrong of course), she
returns to the channel selection state.  If, on the other hand, the
user believes that there is no competing user, she stays at the
access state for another transmission period that has $L_t' \ge L_t$
slots.

The detailed specification of ASA now reduces to finding appropriate
durations of sensing or transmission periods and constructing a
detector for channel occupancy.

\subsection{Time Structure of ASA and Detection Period}\label{Sec:detectionlength}
\begin{figure}[b]
\centering
\begin{psfrags}
\psfrag{P}[c]{\footnotesize Detection period}
\psfrag{C}[c]{\footnotesize sensing}
\psfrag{E}[c]{\footnotesize access}
\includegraphics[width=3.4in]{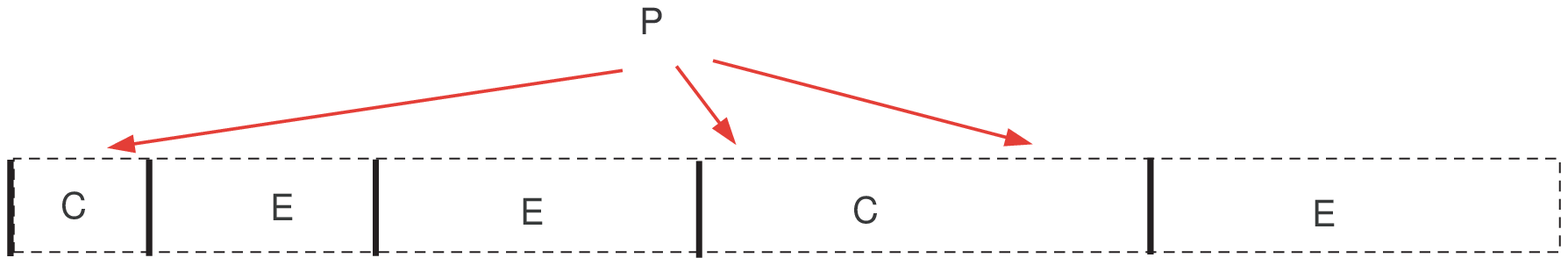}
\end{psfrags}
\caption{Illustration of detection period and increasing detection period length}
\label{fig:pp}
\end{figure}

ASA alternates between sensing and access periods, punctuated by
detection actions.  This structure is illustrated in
Fig.~\ref{fig:pp} where we refer to the time after a detection and
before completing the next detection as a {\em detection period}
during which the user collects either feedback samples (if in the
access state) or observation samples (if in the sensing states)
before a test is performed at the end of the detection period.  The
length of the $k$th detection period is denoted by $L_k$.

A key idea of ASA is to let $L_k$ be a monotonically increasing
function of $k$. Indeed, one can show that if $L_k$ does not grow,
ASA does not achieve the performance offered by the optimal
centralized scheme, due to the non diminishing detection error
probability. Here we choose the form of $L_k$ to be a linear
progression given by
\begin{equation} \label{eq:Lk}
L_{k+1}=L_k + C
\end{equation}
for some integer $C>0$.

The significance of monotonically increasing the detection period is
twofold.  First, with increasing $L_k$, detection accuracy improves.
We show later in Section \ref{sec:schemeanalysis} that the detector
used in ASA has error probabilities decaying exponentially with
respect to $L_k$.

Second, the increasing $L_k$ provides a level of stability to the
policy. A user who finds the correct channel tends to stay there
until completion; it is unlikely a new user can bump a settled user
off its track as time goes.

\subsection{Channel Occupancy Test}
We say that the channel is available in a slot if the channel is on
and no one transmits in that slot.  We present here a detector that
takes channel availability samples and tests channel occupancy, where the two hypotheses are
$\mathcal{H}_0$ channel unoccupied, and $\mathcal{H}_1$ channel occupied by some user.

Note that because the feedbacks from transmission and the outcomes of
sensing give the same information, the detector used at both the
sensing and access states is identical. In both cases, for a sensing
or access period of $L$ slots, the user obtains a sequence of binary
outcomes $\{a, \bar{a}\}^L$ with $a$ indicating that channel is
available.

Let $L_a$ be the number of slots that channel $i$ is available. The
channel occupancy test is a threshold test on the sample mean of the
average availability, \ie
\begin{equation} \label{eq:det}
\frac{L_a}{L}~~ \begin{array}{c}\mathcal{H}_0\mbox{\small: unoccupied}\\ \gtrless\\
\mathcal{H}_1\mbox{\small: occupied}\\ \end{array}~~\eta_i-\epsilon
\end{equation}
where $\epsilon>0$ is a small constant, which lower bounds half of the minimum throughput target ever used by any user in the system (\textit{i.e.}, $\epsilon<r_{\min}/2$)\footnote{$r_{\min}$ is assumed to be a system parameter, which indicates the minimum targeted throughput of interest to the users}.

It is not difficult to see intuitively that, if there is a
persistent user occupying channel $i$, the above detector detects
correctly with high probability.  On the other hand, if the channel
is unoccupied, the probability of mistakenly detecting the channel
as occupied should decay with $L$.

When the underlying channel state processes are alternating renewal
processes assumed in Section \ref{sec:model}, we claim the
following:

\begin{lemma}\label{lemma:detection}
The error probabilities of the channel occupancy detector given in
(\ref{eq:det})
 decay exponentially with $L$.
\end{lemma}

The proof of the above lemma is given in the Appendix. Note that the
length $L_k$ of the $k$th detection period increases linearly with $k$, the
above lemma also implies that detection error probabilities also
decay exponentially with the detection period index $k$.

%% file: Analysis.tex
We present in this section the main results and show that ASA
achieves finite expected regret compared with the optimal
centralized scheme with fixed channel allocation.

Define regret $\mathcal{R}_n$ as the difference in the total number
of successful transmissions (summed over all users) between the
centrally coordinated scenario with pre-determined channel
assignment and the distributed multiaccess scheme in the first $n$
detection periods.  With this we state our main result on expected
regret.

\begin{theorem}\label{thm:regret}
Let $\overline{\mathscr{R}}$ be the maximum throughput region
achievable by a central controller for a $K$ user multiaccess system
involving $N$ independent alternating renewal on-off channels. Then
the expected regret for ASA policy converges to a finite value, \ie
\[
\mbbE(\mathcal{R}_n) \rightarrow d~~\mbox{as $n \rightarrow \infty$}
\]
Consequently, the throughput region achievable by the distributed
learning and access policy ASA $\mathscr{R}_{\mbox{\tiny
ASA}}=\overline{\mathscr{R}}$.
\end{theorem}

\vspace{1em}
The proof of Theorem~\ref{thm:regret} is given in full
in the Appendix.  Here we present a sketch that outlines main ideas
behind the proof.

When all of the $K$ users are in the access state in separate
channels, no expected regret is incurred. Therefore the expected
regret is solely incurred in the undesirable configuration in which
there are still some users not in access state or not in a separate
channel. To investigate the undesirable event, the first ingredient
we need is the exponential decay of $P_{i,e}$, the probability that
there are still some users not in access mode or not in a separate
channel in the $i$th detection period, with respect to the detection
period index $i$. This quantifies the probability of the undesirable
event over the evolution of the policy, and is given in Lemma
\ref{lemma:detect}. Lemma \ref{lemma:Pe} serves as a stepping stone
to Lemma \ref{lemma:detect}.

To capture the expected regret in the first $n$ detection periods,
Lemma \ref{lemma:regretbound} provides an upper bound
(\ref{eqn:upperb}) for the expected regret
$\mathbb{E}\mathcal{R}_n$, which involves $P_{i,e}$ (decreasing with
$i$) and $L_i$ (increasing with $i$). By Lemma \ref{lemma:detect}
$P_{i,e}$ decays exponentially fast, while according to the policy
design, the detection period length $L_i$ only increases linearly.

The fast decay of $P_{i,e}$ and relatively slow growth of $L_i$
guarantees that the upper bound (\ref{eqn:upperb}) sums to a finite
value. This completes the proof of Theorem \ref{thm:regret}. In the
Appendix, we list the required lemmas (Lemma \ref{lemma:regretbound}
to \ref{lemma:detect}) and describe the procedure to prove them.

Provided the finite expected regret result, it is relatively
straight forward to show that the ASA policy achieves identical
throughput region with the centralized scheme with fixed channel allocation, \textit{i.e.},
$\mathscr{R}_{\mbox{\tiny ASA}}=\overline{\mathscr{R}}$, by dividing
the time horizon and taking the limit.

%% file: Appendix.tex
\subsection{Proof of Lemma 1}


To show the exponential decay of the detection error probabilities, namely the false alarm probability $\mathbb{P}(\frac{L_a}{L}<\eta-\epsilon\mid\mathcal{H}_0\mbox{\small: unoccupied})$ and
the miss detection probability $\mathbb{P}(\frac{L_a}{L}>\eta-\epsilon\mid\mathcal{H}_1\mbox{\small: occupied})$, we first notice the mean of the detection statistic $\frac{L_a}{L}$ under both $\mathcal{H}_0$ and $\mathcal{H}_1$,
\begin{equation}
\mathbb{E}_0\frac{L_a}{L}=\eta,
\end{equation}
and
\begin{equation}
\mathbb{E}_1\frac{L_a}{L}=\eta-r,
\end{equation}
where $r$ is the targeted throughput of the user occupying the channel.

Since the channel process is alternating renewal between on and off states, we define an on-off renewal period length to be the total length of an on period and an off period, \textit{i.e.}, an on-off renewal period length $X=X_{\mbox{\tiny on}}+X_{\mbox{\tiny off}}$, where $X_{\mbox{\tiny on}}$ and $X_{\mbox{\tiny off}}$ are
distributed as $F^{\mbox{\tiny on}}$ and $F^{\mbox{\tiny off}}$. The partial sum process $S_n$ is defined as
\begin{equation}
S_n=\sum_{i=1}^nX_i,
\end{equation}
and the associated counting process $N_t$,
\begin{equation}
N_t=\max\{i:S_i\leq t\}.
\end{equation}

Specifically, we will upper and lower bound the detection statistic $L_a/L$, and then prove that both the upper and lower bounds have exponentially decaying probability to deviate from their identical mean, thus showing the detection statistic $L_a/L$ must have exponentially decaying tail probability on both ends away from its mean ($\eta$ for $\mathcal{H}_0$ and $\eta(1-q)=\eta-r$ for $\mathcal{H}_1$, where $q$ is the transmission probability of the user occupying the channel, and $r$ is the targeted throughput). With our choice of $\epsilon$, the detection threshold $\eta-\epsilon$ lies in between the two means under $\mathcal{H}_0$ and $\mathcal{H}_1$.

\begin{proof}
We first upper and lower bound the detection statistic $L_a/L$ in Eq. (\ref{eqn:uplow}).
\begin{equation}\label{eqn:uplow}
\frac{\sum_{i=1}^{N_t}A_i}{t}\leq L_a(t)/t\leq\frac{\sum_{i=1}^{N_t+1}A_i}{t},
\end{equation}
where $A_i$ is the number of available slots (``on" channel state and no other user transmitting) experienced by the user in the $i$th on period for the channel.

We have to show that $L_a(t)/t$ converges to its mean with exponentially decaying tail probability under both hypotheses. This can be done if we can show the leftmost and rightmost sides in Eq. (\ref{eqn:uplow}) converge exponentially fast to their expected value, respectively.

We will treat the rightmost side of Eq. (\ref{eqn:uplow}), and the procedure is similar for the leftmost side.
Specifically, rewrite the rightmost side of Eq. (\ref{eqn:uplow})
\begin{equation}\label{eqn:uplow2}
\frac{\sum_{i=1}^{N_t+1}A_i}{t}=\frac{\sum_{i=1}^{N_t+1}A_i}{N_t+1}\frac{N_t+1}{N_t}\frac{N_t}{t}.
\end{equation}

We cite from Theorem 3.3.2 in \cite{Resnick:ASP} standard almost sure convergence result (\ref{eqn:rt1}) in renewal theory
\begin{equation}\label{eqn:rt1}
\lim_{n\to\infty}\frac{N_t}{t}=\frac{1}{\mu^{\mbox{\tiny on}}+\mu^{\mbox{\tiny off}}},
\end{equation}
and Eq. (\ref{eqn:rt2}) follows directly from Eq. (\ref{eqn:rt1})
\begin{equation}\label{eqn:rt2}
\lim_{n\to\infty}N_t=\infty.
\end{equation}

With Eq. (\ref{eqn:rt1}) and (\ref{eqn:rt2}) it is easy to show that as $t$ approaches infinity, the three terms in Eq. (\ref{eqn:uplow}) $\frac{\sum_{i=1}^{N_t+1}A_i}{N_t+1}$, $\frac{N_t+1}{N_t}$ and $\frac{N_t}{t}$ converge almost surely to their expected values $\mu^{\mbox{\tiny on}}$ under $\mathcal{H}_0$ and $\mu^{\mbox{\tiny on}}(1-q)$ under $\mathcal{H}_1$, 1, $\frac{1}{\mu^{\mbox{\tiny on}}+\mu^{\mbox{\tiny off}}}$, respectively.

Now we turn to the claim that the three terms in Eq. (\ref{eqn:uplow}) converge exponentially fast to their expected values, respectively.

Due to the nature of the alternating renewal channel process and the probabilistic transmissions, the claim for the term $\frac{\sum_{i=1}^{N_t+1}A_i}{N_t+1}$ follows from the standard large deviation result of i.i.d. sum, and Eq. (\ref{eqn:rt2}).

The claim for the terms $\frac{N_t+1}{N_t}$ and $\frac{N_t}{t}$ follows from the assumption in Section \ref{sec:model} that
$\mathbb{E}e^{\theta U^{\mbox{\tiny on}}}$($\mathbb{E}e^{\theta U^{\mbox{\tiny off}}}$) is well defined for some $\theta>0$, where $U^{\mbox{\tiny on}}$($U^{\mbox{\tiny off}}$) is
distributed as $F_i^{\mbox{\tiny on}}$($F_i^{\mbox{\tiny off}}$).
The assumption guarantees that $\mathbb{E}e^{\theta X}$ is well defined for $\theta$, where $X=U^{\mbox{\tiny on}}+U^{\mbox{\tiny off}}$ is the on-off renewal period length.

Before we proceed, cite from Eq. (3.3.1) in \cite{Resnick:ASP} the standard result (\ref{eqn:rt3}) in renewal theory relating the partial sum process with the counting process
\begin{equation}\label{eqn:rt3}
\mathbb{P}(N_t\leq n)=\mathbb{P}(S_n>t).
\end{equation}

Using Eq. (\ref{eqn:rt3}), one has
\begin{eqnarray*}
\mathbb{P}(\frac{N_t+1}{N_t}>1+\epsilon)&=&\mathbb{P}(N_t<\frac{1}{\epsilon})\leq\mathbb{P}(N_t\leq\lceil\frac{1}{\epsilon}\rceil) \\
&=& \mathbb{P}(S_{\lceil\frac{1}{\epsilon}\rceil}>t)= \mathbb{P}(e^{\theta S_{\lceil\frac{1}{\epsilon}\rceil}}>e^{\theta t}) \\
&\leq& \frac{\mathbb{E}e^{\theta S_{\lceil\frac{1}{\epsilon}\rceil}}}{e^{\theta t}}=\frac{\mathbb{E}e^{\theta \sum_{i=1}^{\lceil\frac{1}{\epsilon}\rceil}X_i}}{e^{\theta t}}=\frac{(\mathbb{E}e^{\theta X_i})^{\lceil\frac{1}{\epsilon}\rceil}}{e^{\theta t}}
\end{eqnarray*}
Therefore $\mathbb{P}(\frac{N_t+1}{N_t}>1+\epsilon)$ decays exponentially with respect to $t$.

The claim for the term $\frac{N_t}{t}$ can be treated similarly. We can show the other side of the probability inequalities for $\mathbb{P}(\frac{N_t+1}{N_t}<1-\epsilon)$ and $\mathbb{P}(\frac{N_t}{t}<\frac{1}{\mu^{\mbox{\tiny on}}+\mu^{\mbox{\tiny off}}}-\epsilon)$ in the same way. Therefore we have established the exponential decay of the tail probability of the three terms in Eq. (\ref{eqn:uplow2}). This leads to the exponential decay of the tail probability of the detection statistic $L_a/L$ in the length of the detection period length $t$ (also $L$).

By the threshold structure of the detection, we conclude that the miss detection and false alarm probabilities decay exponentially
with respect to the detection period length $L$.
\end{proof}

\subsection{Lemmas for Theorem \ref{thm:regret}}\label{sec:proof}

We would like to show that the expected regret $\mathbb{E}\mathcal{R}_n$ between the ASA policy and that with central coordination in the first $n$ detection periods converges to a finite value to prove Theorem \ref{thm:regret}. As discussed earlier, the regret $\mathcal{R}_n$ will be small if the fraction of time the users spend in the access mode with orthogonalized channel occupancy is high. Indeed, this relationship is formalized in Lemma \ref{lemma:regretbound}, showing that the expected regret $\mathbb{E}\mathcal{R}_n$ is always upper bounded by (\ref{eqn:upperb}).
\begin{lemma}\label{lemma:regretbound}
\begin{equation}\label{eqn:upperb}
\mathbb{E}\mathcal{R}_n\leq\sum_{i=1}^nNL_iP_{i,e}.
\end{equation}
\end{lemma}

\begin{proof}
We break down the regret according to the detection periods. In the $i$th detection period, if the users are orthogonal and all in access mode, then there is no expected regret incurred. Otherwise, the regret in the $i$th detection period
can at most be as large as the total number of slots contained in the $N$ channels in this detection period, which is exactly the number $NL_i$. Therefore the regret incurred in the $i$th detection period is at most $NL_iP_{i,e}$. Summing over all detection periods from 1 to $n$ yields the desired Eq. (\ref{eqn:upperb}).
\end{proof}

From Eq. (\ref{eqn:upperb}) the expected regret will be finite if $P_{i,e}$, the probability that there are still some users not in access mode or not in a separate channel, decays fast enough compared with the growth of $L_i$.

The factors that drive the decay rate of $P_{i,e}$ include the decay rate of the detection error (how accurate is the inference) and the transition rule's ability to adjust and separate when collisions happen (whether the distributed transition rule is indeed leading users to separate gradually). Lemma \ref{lemma:Pe} shows the quantitative relationship between these two drivers and $P_{i,e}$.

\begin{lemma}\label{lemma:Pe}
The following recursion in the detection period index $i$ holds for $P_{i,e}$
\begin{equation}\label{eqn:PePfm}
P_{i+3,e}\leq3NP_{i,f,m}+P_{i,e}(1-\frac{1}{2^K}\prod_{k=1}^K\frac{N_k-k+1}{N_k}),
\end{equation}
where $P_{i,f,m}$ is the sum of the miss detection probability and the false alarm probability with detection period length $L_i$, $N_k$ is the number of qualified channels for user $k$ (channel ``on" long term fraction no less than the targeted throughput $r_k$), and $N_k$ is ranked in increasing order with $k$.
\end{lemma}

\begin{proof}
The proof of Lemma \ref{lemma:Pe} involves three parts. The first part shows that it is always \textit{possible} that the configuration of the $N$ users will be corrected in at most three detection periods, if the configuration of the current detection period is not orthogonal accessing, and the detection results within the three detection periods are all correct. This part verifies that the transition rule adopted is indeed capable of adjusting and separating the users when collisions happen.

The second part shows that provided that the associated inferences are all correct, the probability of the correction within three detection periods is always larger than $\frac{1}{2^K}\prod_{k=1}^K\frac{N_k-k+1}{N_k}$, where $N_k$ is the number of qualified channels for user $k$, and $N_k$ increases with $k$.
This part verifies that the random channel selection in the transition rule is making strict progress gradually.

The third part verifies Eq. (\ref{eqn:PePfm}) by analyzing the events of detection error and configuration error.

We start by showing the first part by enumerating the possible undesirable configurations.
\begin{enumerate}
\item Several users are still observing separate vacant channels in the $i$th detection period.

 In this case, the user will correctly identify the opportunity and in the $(i+1)$th detection period transition to the access mode.
\item Several users are observing the same vacant channel in the $i$th detection period.

In this case, the user will correctly identify the vacancy and in the $(i+1)$th detection period transition to the access mode. However, this will lead to multiple users transmitting in one channel, which will be correctly detected. Therefore in the $(i+2)$th detection period the users will randomly select channels. With positive (lucky) probability, the selected channels will be vacant and orthogonal, and in the $(i+3)$th detection period the users will transition to the access mode.
\item Several users are accessing the same channel in the $i$th detection period.

In this case, the user will correctly identify the collision among users and in the $(i+1)$th detection period transition to the sensing mode. Therefore in the $(i+1)$th detection period the users will randomly select channels. Still with positive (lucky) probability, the selected channels will be vacant and orthogonal, and in the $(i+2)$th detection period the users will transition to the access mode.
\item Some user is still observing an occupied channel in the $i$th detection period.

In this case, there are two scenarios to analyze: 1) there is another vacant channel qualified for the user, 2) there is currently no vacant channel qualified for the user, \textit{i.e.}, all qualified channels for the user is currently occupied. For scenario 1, the user will correctly identify the fact that the channel is occupied and flip a coin with tail outcome ($1/2$ probability), and in the $(i+1)$th detection period randomly select another qualified channel. With positive (lucky) probability, the selected channel will be vacant, and in the $(i+2)$th detection period the user will transition to the access mode. For scenario 2, the user will correctly identify the fact that the channel is occupied and flip a coin with head outcome ($1/2$ probability), and in the $(i+1)$th detection period start accessing the channel. This will lead to collision in this channel, and in the $(i+2)$th detection period the users will evacuate from the channel and randomly select channels to sense. At this time, with positive (lucky) probability, the selected channels will be vacant and qualified for the users involved in the collision, and in the $(i+3)$th detection period the users will transition to the access mode in separate channels.
\end{enumerate}
Thus we have shown that if the channel configuration of the current detection period is not orthogonal accessing, there is always possibility that the configuration is corrected in at most three steps, provided the detection is correct in the $i$th, $(i+1)$th and $(i+2)$th detection periods.

The second part can be shown by inspecting the required lucky probability. The distributed channel selection by the users incorporates randomness, as well as randomness from the fair coin flipping in the sensing state, and with probability at least $\frac{1}{2^K}\prod_{k=1}^K\frac{N_k-k+1}{N_k}$ the channels can be orthogonalized by the distributed random channel selection, where the factor $\frac{1}{2^K}$ stands for the appropriate fair coin flip probability, and the factor $\prod_{k=1}^K\frac{N_k-k+1}{N_k}$ accounts for the uniform channel selection among qualified channels for each user.

The third part involves analyzing events and algebra. Specifically,
\begin{eqnarray*}
\mathbb{P}_{i+3,e}&\leq&\mathbb{P}(\mathscr{A}_1)+\mathbb{P}(\mathscr{A}_2) \\
&\leq&3NP_{i,f,m}+\mathbb{P}(\mathscr{A}_3) \\
&\leq&3NP_{i,f,m}+P_{i,e}(1-\frac{1}{2^K}\prod_{k=1}^K\frac{N_k-k+1}{N_k}),
\end{eqnarray*}
where $\mathscr{A}_1$ corresponds to the event that at least one user makes a detection error (either miss detection or false alarm in either mode) in the $i$th, $(i+1)$th or $(i+2)$th detection periods,
$\mathscr{A}_2$ corresponds to the event that all detections made by all users are correct in the $i$th, $(i+1)$th and $(i+2)$th detection periods, and the configuration in the $(i+3)$th detection period is still undesirable, $\mathscr{A}_3$ corresponds to the event that the configuration in the $i$th detection period is undesirable, and the random channel selection is not able to separate the users (unlucky), and $P_{i,f,m}$ is the probability that either miss detection or false alarm occurs in one user with detection period length $L_i$.

Specifically, the union bound and the fact that $P_{i+1,f,m}\leq P_{i,f,m}$ yields $\mathbb{P}(\mathscr{A}_1)\leq3NP_{i,f,m}$. Event $\mathscr{A}_2$ is a subset of event $\mathscr{A}_3$, since provided that all detections made by all users are correct in the $i$th, $(i+1)$th and $(i+2)$th detection periods, if either the configuration in the $i$th detection period is desirable, or the random channel selection is able to separate the users (lucky), then the configuration in the $(i+3)$th detection period has to be desirable. This will yield $\mathbb{P}(\mathscr{A}_2)\leq\mathbb{P}(\mathscr{A}_3)$. Finally, the probability that the random channel selection is able to separate the users (lucky) is lower bounded by $\frac{1}{2^K}\prod_{k=1}^K\frac{N_k-k+1}{N_k}$. Therefore the probability that the random channel selection is not able to separate the users (unlucky) is upper bounded by $1-\frac{1}{2^K}\prod_{k=1}^K\frac{N_k-k+1}{N_k}$, and $\mathbb{P}(\mathscr{A}_3)\leq P_{i,e}(1-\frac{1}{2^K}\prod_{k=1}^K\frac{N_k-k+1}{N_k})$.
\end{proof}

With Eq. (\ref{eqn:PePfm}), we are in position to drive the exponential decay of $P_{i,e}$, as established in Lemma \ref{lemma:detect}.
\begin{lemma}\label{lemma:detect}
The probability the system is not in ``good configuration'' in the $i$th detection period, $P_{i,e}$, decays exponentially in the detection period index $i$.
\end{lemma}

\begin{proof}
Write
\begin{equation*}
e^{-3\theta_e}=1-\frac{1}{2^K}\prod_{k=1}^K\frac{N_k-k+1}{N_k}.
\end{equation*}
We have the exponential decay of miss detection probability $P_{i,m}$ and false alarm probability $P_{i,f}$
with respect to the detection period length $L_i$ (linear in the detection period index $i$), which further indicates the exponential decay of the quantity $P_{i,f,m}$ with respect to the detection period index $i$.

Therefore there exists integer $I$, such that for all $i>I$, it holds that
\begin{equation*}
P_{i,f,m}\leq c_{f,m}e^{-\theta_{f,m}i}\leq c_{f,m}e^{-\min\{\theta_{f,m},\frac{\theta_e}{2}\}i}.
\end{equation*}

Write
\begin{equation*}
\alpha=\frac{3Nc_{f,m}}{e^{-3\min\{\theta_{f,m},\frac{\theta_e}{2}\}}-e^{-3\theta_e}},
\end{equation*}
where $\alpha>0$.
Manipulating the recursion equation (\ref{eqn:PePfm}) yields
\begin{eqnarray*}
& & (P_{i+3,e}-\alpha e^{-(i+3)\min\{\theta_{f,m},\frac{\theta_e}{2}\}})-e^{-3\theta_e}(P_{i,e}-\alpha e^{-i\min\{\theta_{f,m},\frac{\theta_e}{2}\}}) \\
&=& P_{i+3,e}-(3Nc_{f,m}e^{-\min\{\theta_{f,m},\frac{\theta_e}{2}\}i}+P_{i,e}e^{-3\theta_e}) \\
&\leq& P_{i+3,e}-(3NP_{i,f,m}+P_{i,e}e^{-3\theta_e})\leq0,
\end{eqnarray*}
which leads to
\begin{equation*}
P_{i+3,e}-\alpha e^{-(i+3)\min\{\theta_{f,m},\frac{\theta_e}{2}\}}\leq e^{-3\theta_e}(P_{i,e}-\alpha e^{-i\min\{\theta_{f,m},\frac{\theta_e}{2}\}}).
\end{equation*}
Therefore one concludes that the exponential decay in $P_{i,e}$ in the detection period index $i$ is guaranteed.
\end{proof} 